\documentclass[10pt]{article}
\usepackage{graphicx}

\def\R{{\rm I\! R}}
\def\N{{\rm I\! N}}

\def \pTW*{\partial_{W^*} T}
\def \pTnW*{\partial^{(n)}_{W^*} T}

\def \dst{\displaystyle}

\def \div{{\rm div\, }}
\def \det{{\rm det\, }}
\def \ra{\rightarrow }

\newtheorem{theorem}{Theorem}
\newtheorem{corollary}{Corollary}

\title{Measure-preserving symmetries and reversibilities of ordinary differential systems} 

\author{Marco Sabatini 
\footnote{Dipartimento di Matematica, Univ. di Trento, I-38123 Povo (TN) - Italy; email: marco.sabatini@unitn.it. This paper has been supported by GNAMPA, Gruppo Nazionale per l'Analisi Matematica, la Probabilit\`a e le loro Applicazioni.}
}
\begin{document}
\maketitle
\begin{abstract}  We prove that measure-preserving symmetries of an $n$-dimensional differential system preserve its divergence and the divergence derivatives along the solutions. Also, we prove that measure-preserving reversibilities preserve odd-order divergence derivatives along the solutions, and that even-order derivatives are multiplied by $-1$.  We apply such results to find all the area-preserving symmetries and reversibilities of planar Lotka-Volterra and Li\'enard systems. 

  \end{abstract}

%%%%%%%%%%%%%%

%%%%%%%%%%%%%%

\section{Introduction}

Let us consider a  differential system
\begin{equation} \label{sysn}
\dot z = F(z),
\end{equation}  
where $F(z)=(F_1(z),\dots, F_n(z)) \in C^\infty (\Omega,\R^n)$,  $\Omega \subset \R^n$ open and connected, $z=(z_1, \dots, z_n) \in \Omega$. We denote by $\phi(t,z)$ the local flow defined by (\ref{sysn}) in $\Omega$.  
Differential systems are often studied with the aid of suitable transformations, either to put them into a more convenient form, or to reduce their study to a proper portion of the  space. In some cases the second option can lead to discover dynamic properties of the system otherwise difficult to prove. This is the case of some  symmetries, as for instance  mirror symmetries, that can show the existence of periodic solutions in presence of rotating orbits. For planar systems this is the easiest way to prove the existence of cycles in absence of known first integrals \cite{Col,GM,R,vW}.  Several papers were devoted to finding conditions for the existence of mirror symmetries, in particular for polynomial systems \cite{ACGG,BL,Col,GM,R,S2,vW,Z}. Such symmetries are usually called reversibilities, and are characterized by the following flow property,
\begin{equation} \label{flowrev}
\sigma(\phi(t,z)) = \phi(-t,\sigma(z)),
\end{equation} 
where $\sigma$ is the mirror symmetry. The existence of such a reversibility can be proved without knowing the flow,  by showing that 
\begin{equation} \label{Jrev}
  F(\sigma(z)) = -J_\sigma(z) \cdot F(z),
\end{equation} 
where $J_\sigma$ is the Jacobian matrix of $\sigma$.

The argument showing that every rotating orbit of a planar system that meets the symmetry axis is a cycle applies without changes to mirror-like nonlinear symmetries satisying (\ref{flowrev}). This provides  a more general tool for proving local integrability, even if finding mirror-like nonlinear symmetries is much more difficult than finding mirror ones. Some progress has been made recently in relation to measure-preserving symmetries  \cite{S1}, proving that the symmetry hypersurface is contained in the  null-divergence set  of (\ref{sysn}). 

\medskip

A second class of symmetries, characterized by a flow commutativity property,
\begin{equation} \label{flowsymm}
\sigma(\phi(t,z)) = \phi(t,\sigma(z)),
\end{equation} 
 is also useful in studying differential systems \cite{S0}. 
This is the case of flows generated by point symmetric vector fields. Even the existence of such  symmetries can be proved without knowing the flow,   by showing that 
\begin{equation} \label{Jsymm}
  F(\sigma(z)) = J_\sigma(z) \cdot F(z).
\end{equation} 

\medskip

In this paper we are concerned with measure-preserving nonlinear symmetries of both types. 
Measure-preserving transformations occur in several fields of mathematics. For instance, every non-singular linear transformation can be normalized to a measure-preserving linear transformation. Another class of measure-preserving transformations is given by triangular maps,
$$
(x_1, \dots, x_n) \mapsto \Big(x_1,x_2+f_1(x_1), x_3+f_2(x_1,x_2), \dots,x_n + f_{n-1}(x_1, \dots, x_{n-1}) \Big) .
$$
The polynomial maps considered by the celebrated Jacobian Conjecture \cite{ BCW} are measure-preserving, too. Nonlinear rotations, defined by
$$
  \matrix{\sigma(x,y) = \Big(x \cos \alpha \left( \rho \right) - y \sin  \alpha \left( \rho \right) ,  x \sin \alpha \left( \rho \right) + y \cos  \alpha \left( \rho \right) \Big) ,}  
$$ 
where $\displaystyle{  \rho = \sqrt{x^2 + y^2}   }  $ and  $\alpha(\rho)$ is a  smooth scalar function, provide another class of planar measure-preserving transformations.

\medskip

For the sake of simplicity, in the following we call symmetries only the transformations satisfying (\ref{flowsymm}). In other words, in this paper terminology, symmetries are associated to commutativity (\ref{flowsymm}), reversibilities to anti-commutativity (\ref{flowrev}). 

\medskip

Denoting by $D(z) \equiv \div F(z)$ the divergence $F$ at $z$, for all $j\in \N$  we denote by
$$ D^{(j)}(z) =  \frac{d^j D(\phi(t,z))}{d t^j}\bigg|_{t=0}
$$
the $j$-th derivative of the divergence along the local flow. 
In this paper we prove that for a measure-preserving reversibility $\sigma$ one has, for all $j\in \N$, 
$$
D^{(j)}  (\sigma(z)) =  (-1)^{j+1} D^{(j)}(z)  .
$$
As a consequence we extend the main result proved in \cite{S1}, showing that at a $\sigma$-fixed point $z$ one has $D^{( 2j)}  (z) =0$.

Similarly,  we prove that  for a measure-preserving symmetry $\sigma$ one has, for all $j\in \N$, 
$$
D^{(j)}  (\sigma(z)) =  D^{(j)}(z)  .
$$
Such results allow us to give necessary conditions for the existence of symmetries and reversibilities. Also,  we give the implicit form of all the measure-preserving symmetries and reversibilities a system can have. 

\medskip

The above results are applied to classify  the area-preserving symmetries and reversibilities of planar Lotka-Volterra and Li\'enard systems. As for Li\'enard systems, we assume standard hypotheses, usually considered in relation to the study of rotation points.

\medskip

 For the sake of simplicity we assume all the functions and vector fields to be of class $C^\infty$, but  all  the results presented here can be proved under weaker regularity hypotheses.

\section{Symmetries and reversibilities in $\R^n$}

Let $\sigma: \Omega \rightarrow \Omega$ be a homeomorhism such that $\sigma^2(z) = \sigma(\sigma(z)) =z$. Such maps are usually called {\it involutions}. If $A \subset \Omega$, we say that $A$ is {\it $\sigma$-invariant} if $\sigma(A) = A$. 
A point $z \in \Omega$ is said to be a {\it $\sigma$-fixed point} if $\sigma(z)=z$.  Elementary examples of involutions are mirror symmetries and point symmetries. Mirror symmetries have a symmetry hypersurface consisting of the symmetry fixed points. Point symmetries have a symmetry center consisting of the symmetry fixed point. 
Denoting by $\mu$ the $n$-dimensional Lebesgue measure, we say that $\sigma$ is a measure-preserving involution if, for every measurable set $A$, one has
$$
\mu\big(\sigma(A)\big) = \mu (A).
$$

We say that the system (\ref{sysn}) is {\it $\sigma$-symmetric} if  there exists an involution $\sigma$ such that
$$
\sigma(\phi(t,z)) = \phi(t,\sigma(z)),
$$
for all $t$ and $z$ such that the above equality makes sense.
Similarly, we say that the system (\ref{sysn}) is {\it $\sigma$-reversible} if  there exists an involution $\sigma$ such that
$$
\sigma(\phi(t,z)) = \phi(-t,\sigma(z)),
$$
Such  properties can be revealed by verifying an equality involving the vector field. In fact,  if $\sigma \in C^\infty(\Omega,\Omega)$, then the $\sigma$-symmetry  of (\ref{sysn}) can be checked by verifying the relationship
$$
 F(\sigma(z)) =  J_\sigma(z) \cdot F(z) ,
$$
where $J_\sigma(z)$ is the Jacobian matrix of $\sigma$ at $z$. Similarly, the  $\sigma$-reversibility is equivalent to
$$
 F(\sigma(z)) =  -J_\sigma(z) \cdot F(z) ,
$$

\medskip

 If $\sigma$ is a  symmetry or a reversibility for the local flow $\varphi$, then it takes equilibrium points into equilibrium points and cycles into cycles.
 In general,  the fixed sets of $\varphi$ and $\sigma$ do not coincide. For instance the system
\begin{equation}  
\left\{  \matrix{\dot x  =    y \hfill  \cr    \dot y =   - \sin x \hfill,}  \right.
\end{equation}
is symmetric with respect to $\sigma(x,y) = (-x,-y)$ which has a single fixed point at the origin. Such a system has infinitely many critical points. On the other hand, if a system has a unique critical point, then it is $\sigma$-fixed with respect to every symmetry or reversibility.

\medskip

Denoting by $D(z)$ the divergence $\div F(z)$ of $F$ at $z$, for all $j\in \N$  we denote by
$$ D^{(j)}(z) =  \frac{d^j D(\phi(t,z))}{d t^j}\bigg|_{t=0}
$$
the $j$-th derivative of $D(z)$ along the solutions of (\ref{sysn}). If $j=0, 1, 2$ we  write  $ D(z)$, $\dot D(z)$, $\ddot D(z)$, rather than $D^{(0)}(z) $, $D^{(1)}(z) $, $D^{(2)}(z) $.

\medskip

In next proof, for $\varepsilon > 0$ we denote by $\Sigma_\varepsilon$ be a $(n-1)$-dimensional open hypercube of sidelength $\varepsilon$.

\begin{theorem} \label{SRD} Let $\sigma\in C^\infty(\Omega,\Omega)$ be a measure-preserving involution. Assume (\ref{sysn}) to have isolated critical points.
The following statements hold.
\begin{itemize}
\item[1)] If $\sigma$ is a symmetry of (\ref{sysn}), then for all $j\in \N$, 
$$D^{(j)}  (\sigma(z)) = D^{(j)}(z)   .$$
\item[2)] If $\sigma$ is a reversibility of (\ref{sysn}), then for all $j\in \N$, 
$$D^{(j)}  (\sigma(z)) = (-1)^{j+1} D^{(j)}(z)  . $$
\end{itemize}
 \end{theorem}
{\it Proof.}    $1)$
By continuity, it is sufficient to prove the statement for non-critical points.

Let us proceed  by induction. Let $z\in \Omega$ be such that $F(z) \neq 0$.  For $\varepsilon > 0$, let $\Sigma_\varepsilon$ be transversal to $F(z)$ at $z$. Since $\sigma$ is a diffeomorphism, the curve $\sigma( \Sigma_\varepsilon)$ is transversal to $F(\sigma(z))$ at $\sigma(z)$.
There exists $\varepsilon > 0$ small enough to have $A_\varepsilon:= \phi((-\varepsilon,\varepsilon),\Sigma_\varepsilon) \subset \Omega$ and 
$\sigma(A_\varepsilon)  \subset \Omega$. Similarly, there exists $t^* > 0$ small enough to have 
$\phi(t, A_\varepsilon)  \subset \Omega$ and $\sigma( \phi(t,A_\varepsilon)) = \phi (t,\sigma(A_\varepsilon))  \subset \Omega$ for all $t \in (-t^*,t^*)$.

By Liouville theorem and the integral mean value theorem there exists $\chi_{t \varepsilon }\in \phi(t,A_\varepsilon)$ such that:
$$
\frac{d \mu(\phi(t,A_\varepsilon))}{dt} = \int_{\phi(t,A_\varepsilon)} \!  \div F(z)\ dz = \mu(\phi(t,(A_\varepsilon))\, \div F(\chi_{t \varepsilon}).
$$
Similarly, there exists $\eta_{t \varepsilon }\in \phi (t,\sigma(A_\varepsilon)) = \sigma( \phi(t,A_\varepsilon))$ such that:
$$
 \frac{d \mu(\phi(t,\sigma(A_\varepsilon)))}{dt} = \int_{\phi (t,\sigma(A_\varepsilon))}\!  \div F(z)\ dz = \mu(\phi(t,\sigma(A_\varepsilon)))\, \div F(\eta_{t \varepsilon }).
$$
Since $\sigma$ is a measure-preserving symmetry one has:
$$
  \frac{d \mu(\phi(t,\sigma(A_\varepsilon)))}{dt}\bigg|_{t=0} = \lim_{t \to 0} \frac {\mu(\phi(t,\sigma(A_\varepsilon))) - \mu(\sigma(A_\varepsilon))}{t}   = \lim_{t \to 0} \frac {\mu(\sigma(\phi(t,A_\varepsilon))) - \mu(\sigma(A_\varepsilon))}{t} = 
$$  $$
  =\lim_{t \to 0} \frac {\mu(\phi(t,A_\varepsilon)) - \mu(A_\varepsilon)}{t} =
 \frac{d \mu(\phi(t,A_\varepsilon))}{dt} \bigg|_{t=0}   .
 $$
 Hence one has, for $t = 0$,
 $$
 \mu(A_\varepsilon)\, \div F(\chi_{0 \varepsilon}) =  \frac{d \mu(\phi(t,A_\varepsilon))}{dt} \bigg|_{t=0}  =
    \frac{d \mu(\phi(t,\sigma(A_\varepsilon)))}{dt}\bigg|_{t=0} = 
 \mu(\sigma(A_\varepsilon))\, \div F(\eta_{0 \varepsilon }) . 
 $$
 which implies $\div F(\chi_{0 \varepsilon}) = \div F(\eta_{0 \varepsilon })$, since $ \mu(A_\varepsilon) = \mu(\sigma(A_\varepsilon))$. As $\varepsilon $ tends to $0$,  $\chi_{0 \varepsilon}$ tends to $z$ and $\eta_{0 \varepsilon }$ tends to $\sigma(z)$, so that by the continuity of $\div F(z)$ one has
 $$
\div F(\sigma(z)) = \div F(z).
 $$
 Now assume by induction that $D^{(j)}  (\sigma(z)) = D^{(j)}(z) $. One has
$$
D^{( j+1)}  (\sigma(z)) =  \lim_{t \to 0} \frac {D^{(j)}(\phi(t,\sigma(z))) - D^{(j)}(\sigma(z)) }{t} =
  \lim_{t \to 0} \frac {D^{(j)}(\sigma(\phi(t,z))) - D^{(j)}(\sigma(z)) }{t} =
  $$  $$
=  \lim_{t \to 0} \frac {D^{(j)}(\phi(t,z)) - D^{(j)}(z) }{t} =D^{( j+1)}  (z).
$$ 
 
 $2)$ Analogous to 1). In particular, 
 since $\sigma$ is a measure-preserving reversibility one has:
 $$
  \frac{d \mu(\phi(t,\sigma(A_\varepsilon)))}{dt}\bigg|_{t=0} = \lim_{t \to 0} \frac {\mu(\phi(t,\sigma(A_\varepsilon))) - \mu(\sigma(A_\varepsilon))}{t}   = \lim_{t \to 0} \frac {\mu(\sigma(\phi(-t,A_\varepsilon))) - \mu(\sigma(A_\varepsilon))}{t} = 
$$  $$
 =  - \lim_{t \to 0} \frac {\mu(\sigma(\phi(-t,A_\varepsilon))) - \mu(\sigma(A_\varepsilon))}{-t} = 
- \lim_{t \to 0} \frac {\mu(\phi(-t,A_\varepsilon)) - \mu(A_\varepsilon)}{-t} = 
- \frac{d \mu(\phi(t,A_\varepsilon))}{dt} \bigg|_{t=0}   .
 $$
 This gives, working as in the previous case,
  $$
  \div F(\sigma(z)) = - \div F(z) ,
 $$
that is $D(\sigma(z))= - D(z)$. 

Then,  assuming by induction that $D^{(j)}  (\sigma(z)) = (-1)^{j+1} D^{(j)}(z)  $, one has
$$
D^{( j+1)}  (\sigma(z)) =  \lim_{t \to 0} \frac {D^{(j)}(\phi(t,\sigma(z))) - D^{(j)}(\sigma(z)) }{t} =
  $$  $$
=    \lim_{t \to 0} \frac {D^{(j)}(\sigma(\phi(-t,z))) - D^{(j)}(\sigma(z)) }{t} =
\lim_{t \to 0} \frac {(-1)^{j+1} D^{(j)}(\phi(-t,z)) -(-1)^{j+1}  D^{(j)}(z) }{t} =
$$  $$
=- (-1)^{j+1} \lim_{t \to 0} \frac { D^{(j)}(\phi(-t,z)) - D^{(j)}(z) }{-t} 
= (-1)^{j+2} D^{( j+1)}  (z).
$$ 
\hfill  $\clubsuit$ \\

The theorem \ref{SRD} can be read as a statement regarding invariance properties of $D^{(j)} $ level sets, as in next corollary.

\begin{corollary} \label{SRDcor1} Let $\sigma\in C^\infty(\Omega,\Omega)$ be a measure-preserving involution and (\ref{sysn}) have isolated critical points. 
 The following statements hold.
\begin{itemize}
\item[1)] If $\sigma$ is a symmetry  of (\ref{sysn}), then for all $j\in \N$ every level set 
$  D^{(j)}(z)  = L $ is $\sigma$-invariant.
\item[2)] If $\sigma$ is a reversibility of (\ref{sysn}), then for all $j\in \N$, every level set 
$  D^{(2j+1)}(z)  = L $ is $\sigma$-invariant. Moreover,  for all $j\in \N$, every level set 
$  \left(  D^{(2j)}(z) \right)^2 = L $ is $\sigma$-invariant.\end{itemize}
 In particular, both for symmetries and for reversibilities,  for all $j\in \N$ the set $D^{(j)}(z)  =0 $ is $\sigma$-invariant.
\end{corollary}
{\it Proof.}  Immediate, from theorem \ref{SRD}.
 \hfill  $\clubsuit$ \\

Next corollary extends theorem (1) in \cite{S1}. 

\begin{corollary} \label{SRDcor2}  Let $\sigma\in C^\infty(\Omega,\Omega)$ be a measure-preserving involution and (\ref{sysn}) have isolated critical points.  
 If $\sigma$ is a  reversibility  of (\ref{sysn}), then for all $j\in \N$ the set
$ D^{( 2j)}(z)  =0 $
contains all the fixed points of $\sigma$.
\end{corollary}
{\it Proof.}   If $\sigma(z) = z$, then $ D^{( 2j)}(z)  = D^{( 2j)}(\sigma(z))  = - D^{( 2j)}(z)    $, hence 
 $ D^{( 2j)}(z) =0$.  \hfill  $\clubsuit$ \\
 
 The corollary \ref{SRDcor2} cannot be extended to symmetries. In fact, the system
$$
\dot x =  x ,
$$
is symmetric w. resp. to the involution $\sigma(x,y) = -x $, which has a unique fixed point at the origin, and has $D = 1$.

\medskip

In general the fixed set of $\sigma$ does not coincide with the null divergence set. The system
\begin{equation}  \label{xy2}
\left\{  \matrix{\dot x =  x^2 \cr    \dot y =  y^2 ,} \right.
\end{equation}
is reversible w. resp. to  the involution $\sigma(x,y) = (-x,-y)$, whose fixed set is the origin, properly contained in the null divergence set $x+y=0$.

\bigskip

We say that a strictly increasing map $\pi: \{1,\dots, n\} \ra \N$ is a {\it selection}. Given  a selection, let us define a map $\Delta_\pi$ on $\Omega$ as follows.
\begin{equation} \label{delta}
\Delta_\pi(z) = \left( D^{( \pi(1))}(z), \dots,  D^{( \pi(n))}(z) \right).
\end{equation}
The simplest choice consists in taking $\pi (j)=j-1$, so obtaining
\begin{equation} \label{deltaid}
\Delta_\pi(z) = \left( D(z), \dot D(z), \dots,  D^{( n-1)}(z) \right),
\end{equation}
the map whose components are $D(z)$ and  the first $n-1$ derivatives of $D(z)$ along the solutions of (\ref{sysn}). In this case we omit the subscript $\pi$,  writing only $\Delta(z)$, rather than $\Delta_\pi(z)$.  If $n=2$, one  has
$$
\Delta(z) = \left( D(z), \dot D(z) \right).
$$

\medskip

For every $\pi$ as above, we define the {\it sign matrix $S_\pi $ of $\pi$} as the $n \times n$ diagonal matrix such that
$$
a_{jj} = (-1)^{\pi(j)+1}.
$$
If $\pi(j) = j-1$ one has, for $n = 2,3,4$,
$$
S_\pi =\left(  \matrix{-1   & 0 \cr    0  & 1}  \right), \qquad 
   S_\pi =\left(  \matrix{-1   & 0 & 0 \cr 
   0  & 1 & 0 \cr 0 & 0 & -1}  \right) , \qquad
    S_\pi =\left(  \matrix{-1   & 0 & 0 & 0\cr 
   0  & 1 & 0 & 0\cr 0 & 0 & -1 & 0  \cr 
   0  & 0 & 0 & 1}  \right) .   
$$
The statement of theorem \ref{SRD} can be rephrased as follows.

\begin{corollary} \label{SRDcor3} Let $\sigma\in C^\infty(\Omega,\Omega)$ be a measure-preserving involution. 
The following statements hold.
\begin{itemize}
\item[1)] If $\sigma$ is a symmetry of (\ref{sysn}), then for every selection $\pi$
$$\Delta_\pi(\sigma(z)) = \Delta_\pi(z)  .$$
\item[2)] If $\sigma$ is a reversibility of (\ref{sysn}),then for every selection $\pi$
$$\Delta_\pi(\sigma(z)) =S_\pi \cdot \Delta_\pi(z)  .$$
\end{itemize}
 \end{corollary}

In order to illustrate how the above facts can be applied, let us consider the system
\begin{equation}  \label{cenquad}
\left\{  \matrix{\dot x  =    y + x^2 \cr    \dot y =   -g(x) \hfill  .}  \right.
\end{equation}
One has
$$
\Delta(x,y) = (D, \dot D)  = (2x, 2y + 2x^2).
$$
$\Delta$'s components are independent of $g(x)$. Let us look for conditions for (\ref{cenquad}) to be symmetric or reversible.  As for symmetry, setting $\sigma(x_1,y_1) = (x_2,y_2) $, one has
$$ \left\{
\matrix { 2x_2=D(x_2,y_2) = D(x_1,y_1) = 2x_1 \hfill  \cr 
 2y_2 + 2x_2^2=\dot D(x_2,y_2) = \dot D(x_1,y_1) = 2y_1 + 2x_1^2. }
\right.
$$
From the first equation one has $x_2=x_1$, then from the second one $y_2=y_1$.  This proves that the unique possible area-preserving symmetry (APS) of (\ref{cenquad}) is the identity map $\sigma(x,y) = (x,y) $, regardless of $g(x)$. 

\medskip

On the other hand, if there exists an area-preserving reversibility (APR)  $\sigma$ of  (\ref{cenquad}), then, setting $\sigma(x_1,y_1) = (x_2,y_2) $, one has
$$ \left\{
\matrix { 2x_2=D(x_2,y_2) = -D(x_1,y_1) = -2x_1 \hfill  \cr 
 2y_2 + 2x_2^2=\dot D(x_2,y_2) = \dot D(x_1,y_1) = 2y_1 + 2x_1^2. }
\right.
$$
From the first equation one has $x_2=-x_1$, then from the second one $y_2=y_1$.  Hence the unique candidate to be an APR of (\ref{cenquad}) is the  map $\sigma(x,y) = (-x,y) $. Checking the condition (\ref{Jrev}), one has
 $$
 F(\sigma(x,y)) = F(-x,y) = (y+x^2,-g(-x)), \quad J_\sigma(x,y) \cdot F(x,y) = (-y-x^2,-g(x)).
 $$
The equality $ F(\sigma(z)) =  - J_\sigma(z) \cdot F(z) $ holds identically if and only if $g(x) = - g(-x)$. This proves that if $g(x)$ is odd, then $\sigma(x,y) = (-x,y)$ is a reversibility of the system (\ref{cenquad}).  Moreover, if $g(x)$ is not odd, then such a system has no APRs.

\medskip

Given a selection $\pi$, assume $\Delta_\pi$ to be locally invertible at some  $z\in \Omega$. Let $U_z$ be a domain of invertibility and $\Delta_\pi^z$ be the restriction of $\Delta_\pi$ at $U_z$.

\begin{corollary} \label{charsigma} Let $\sigma\in C^\infty(\Omega,\Omega)$ be a measure-preserving involution and (\ref{sysn}) have isolated critical points. Assume there exist a selection $\pi$ and $z \neq \sigma(z)  $  such that   $\Delta_\pi$ is locally invertible both at $z$ and at $\sigma(z)$.
Then the following statements hold.
\begin{itemize}
\item[1)] If $\sigma$ is a symmetry of (\ref{sysn}), then 
$$\sigma(z) =\left(  \Delta_\pi ^{\sigma(z)}\right)^{-1}  \Big(  \Delta_\pi^z(z) \Big)  .$$
\item[2)] If $\sigma$ is a reversibility of (\ref{sysn}), then 
$$\sigma(z) =\left(  \Delta_\pi ^{\sigma(z)}\right)^{-1} \Big( S_\pi \cdot   \Delta_\pi^z(z) \Big)  .$$
\end{itemize}
 \end{corollary}
 {\it Proof.}   Immediate from Corollary  \ref{SRDcor3}.
 \hfill  $\clubsuit$ \\

Next corollary is concerned with symmetries and reversibilities endowed with a selection whose sign matrix is  the identity.

\begin{corollary} \label{notinv} Let $n\geq 2$, (\ref{sysn}) have isolated critical points and $\sigma\in C^\infty(\Omega,\Omega)$  be a measure-preserving symmetry or reversibility.  Let $\pi$ be a selection such that $\Delta_\pi(\sigma(z))= \Delta_\pi(z) $. If $z_0$ is a $\sigma$-fixed point and $\sigma$ is non-trivial in any neighbourhood of $z_0$, then the map  $\Delta_\pi$ is not locally invertible at $z_0$. \end{corollary}
 {\it Proof.}   
  Assume by absurd $\Delta_\pi$ to be invertible on a neighbourhood $U_0 $ of $z_0$. Since $\sigma$ is non-trivial in , there exists a sequence $w_n\in U_0$ converging to $z_0$, with  $\sigma(w_n) \neq w_n$. By continuity, $\sigma(w_n) $ converges to $\sigma(z_0)=z_0$. Since both sequences, $w_n$ and $\sigma(w_n) $, converge to the same point $z_0$, they definitely belong  to $U_0$. On the other hand, one has $\Delta_\pi(\sigma(w_n))= \Delta_\pi(w_n) $, contradicting $\Delta_\pi$ invertibility on $U_0$.   
 \hfill  $\clubsuit$ \\

The corollary  \ref{notinv} cannot be extended to arbitrary selections of reversibilities. In fact, the system
$$
\left\{  \matrix{\dot x  =    y + x^2 \cr    \dot y =   -x \hfill  }  \right.
$$
is reversible with respect to the involution $\sigma(x,y) = (-x,y)$, which has  infinitely many fixed points.
It has $D=2x$, $\dot D = 2y + 2x^2$, hence $\Delta (x,y) = (2x,2y+2x^2)$, which is globally invertible on $\R^2$.

\section{Planar Lotka-Volterra systems}

In this section we are concerned with symmetries and reversibilities of the following class of systems, 
\begin{equation} \label{sysLV}
\left\{  \matrix{\dot x  =    x(a-by)  \hfill  \cr    \dot y =  y(cx-d) .}  \right.
\end{equation} 
Such systems arise in the study of biomathematical models. In such a case $a, b, c, d$ are positive real numbers. There exist two critical points, the origin $O$ and 
$$
P =\left( \frac  dc, \frac ab \right).
$$
$O$ is a saddle point. The axes are its separatrices, invariant lines bounding  a region covered, if $a, b, c, d >0$, with non-trivial cycles surrounding  the critical point $P$. 
If $a, b, c, d$ are not all positive, we can consider the jacobian matrices of (\ref{sysLV}) at $O$ and $P$,
$$
J(O) = \left(  \matrix{a & 0 \cr   0 & -d\hfill  }  \right),
$$
$$
J(P) = \left(  \matrix{ 0 &  -\frac{bd}{c} \cr   \frac{ac}{b} &  0\hfill  }  \right).
$$
Their determinants are, respectively, $-ad$ and $ad$, hence  $O$ and $P$ are critical points of different type. If a symmetry or reversibility $\sigma$ exists, it takes critical points into critical points of the same type, hence  $O$ and $P$ are both $\sigma$-fixed points.

\medskip

From now on we only consider non-trivial symmetries and reversibilities, without stating it explicitely in the theorems.

\medskip

In next theorem we do not assume any sign conditions on $a, b, c, d$. We just require $b$ and $c$ not to vanish, since if $b=0$ or $c=0$, then the system degenerates into a triangular one,  integrable using scalar linear equation techniques.

\begin{theorem} \label{lvrev} 
The system (\ref{sysLV}), with $bc \neq 0$, has an APR if and only if
$$
a = d.
$$
If such a condition holds, then the unique   APR of (\ref{sysLV}) is 
$$ \sigma(x,y) = \left( \frac {by}c , \frac {cx}b  \right) . $$
\end{theorem}
{\it Proof.}  
Necessity.

One has
$$
D(x,y) =cx - by +a -d.
$$
The set $D=0$ is a line $r$ containing all $\sigma$-fixed points \cite{S1}. For every $a , b , c, d$ one has $P \in r$. On the other hand,  $O \in r$ if and only if $a = d$.

\medskip

Sufficiency. 

Let us  assume $a =d$. Then, one has
$$
\dot D(x,y)  = -2bcxy+a(cx+by).
$$
In order to find a reversibility, we have to solve the following system w. resp. to $(x_2,y_2)$,
$$
\left\{  \matrix{ cx_2 - by_2 = -cx_1 + by_1  \hfill  \cr -2bcx_2y_2+a(cx_2+by_2) = -2bcx_1y_1+a(cx_1+by_1)  }  \right.
$$
There exist two solutions,
$$
   (x_2,y_2) = \left( \frac ac - x_1, \frac ab - y_1 \right) \qquad (x_2,y_2) =\left( \frac {by_1}c , \frac {cx_1}b  \right) .
$$
The map $\dst{ \sigma(x,y) = \left( \frac ac - x, \frac ab - y \right)   }$ takes $O$ into $P$ and vice versa, hence it is not a reversibility of (\ref{sysLV}). On the other hand, since $a = d$, the map $\dst{ \sigma(x,y) = \left( \frac {by}c , \frac {cx}b  \right) }$ is an area-preserving involution that leaves $O$ and $P$ fixed. The reversibility condition is satisfied, since
$$
F(\sigma(x,y)) =  \left(   \frac {aby}c - bxy  , cxy - \frac {acx}b  \right) ,$$
and
$$
- J_\sigma(x,y) \cdot F(x,y) = -
 \left( \matrix{
 0 & \frac bc \hfill  \cr \frac cb & 0
 }  \right) \cdot  \left( \matrix{
  x(a-by)   \hfill  \cr y(cx-a)
 }  \right) = \left(   \frac {aby}c - bxy  , cxy -  \frac {acx}b  \right) .
$$
\hfill  $\clubsuit$ \\

The above theorem  gives reversibilities also for classical planar Lotka-Volterra systems, where all coefficients are positive. On the opposite hand, next theorem gives symmetries only for systems with coefficients of opposite signs.

\begin{theorem} \label{lvsym} 
The system (\ref{sysLV}), with $bc \neq 0$, has an APS if and only if
$$
a + d = 0.
$$
If such a condition holds, then the unique   APR of (\ref{sysLV}) is 
$$ \sigma(x,y) = \left( - \frac {by}c , - \frac {cx}b  \right) . $$
\end{theorem}
{\it Proof.}  
Necessity.

The map $\Delta$ has the following form,
$$
\Delta(x,y) = \big(  cx - by + a - d, -2bcxy+acx + bdy \big) .
$$
Its jacobian matrix is
$$
J_\Delta(x,y) = \left(  \matrix{ c & -b    \cr  c ( a - 2by) & b(d - 2cx ) }  \right).
$$
Its determinant is $\det J_\Delta(x,y) = bc ( a+ d - 2 c x - 2 b y )   $. If an APS  $\sigma$ exists, then $\Delta$  is  not locally invertible at $\sigma$-fixed points. Both $O$ and $P$ are fixed, hence one has
$$
J_\Delta(0,0) = bc(a+d) = 0 , \qquad J_\Delta\left( \frac  dc, \frac ab \right) =    - bc (a+d) = 0           .
$$
In both cases one has $a+d=0$. 
 
 \medskip
 
 Since $d = -a$, in order to find the APS form, we have to solve the following system, equivalent to $\Delta(x_2,y_2) = \Delta(x_1,y_1)$, w. resp. to $(x_2,y_2)$, 
$$
\left\{  \matrix{ cx_2 - by_2 + 2a = cx_1 - by_1 + 2a  \hfill  \cr -2bcx_2y_2+a(cx_2-by_2) = -2bcx_1y_1+a(cx_1-by_1)  }  \right.
$$
From the first equation one has $cx_2 - by_2 = cx_1 - by_1$. A substitution into the second equation leads to  $x_2y_2=x_1y_1$. The resulting system,
$$
\left\{  \matrix{ cx_2 - by_2  = cx_1 - by_1  \hfill  \cr x_2y_2=x_1y_1  }  \right.
$$
has two solutions,
$$
(x_2,y_2) = (x_1,y_1), \qquad (x_2,y_2) = \left(-\frac {b y_1}c,-\frac {cx_1}b \right) .
$$
The first solution corresponds to the trivial APS, the second one to the involution
\begin{equation}  \label{sigmaLVsym}
\sigma(x,y) =  \left(-\frac {b y}c,-\frac {cx}b \right) .
\end{equation}

 \medskip 
 
 Sufficiency. 
 
 Assuming $a+d=0$, let us define $\sigma$ as in (\ref{sigmaLVsym}). Checking the condition (\ref{Jsymm}) one has
 $$
F(\sigma(x,y)) =  \left(   - \frac {aby}c - bxy  , cxy - \frac {acx}b  \right) ,
$$
and
$$
 J_\sigma(x,y) \cdot F(x,y) = 
 \left( \matrix{
 0 & - \frac bc \hfill  \cr - \frac cb & 0
 }  \right) \cdot  \left( \matrix{
  x(a-by)   \hfill  \cr y(cx+a)
 }  \right) = \left(   - \frac {aby}c - bxy  , cxy -  \frac {acx}b  \right) .
$$

\hfill  $\clubsuit$ \\

\section{Li\'enard systems}

In this section we are concerned with reversibilities and symmetries   of the following class of systems, 
\begin{equation} \label{syslie}
\left\{  \matrix{\dot x  =    y  \hfill  \cr    \dot y =   -g(x) - yf(x) , }  \right.
\end{equation} 
equivalent to the Li\'enard equation $\ddot x + f(x) \dot x + g(x) = 0$. One has $\Delta(x,y) = (D(x,y),\dot D(x,y)) = (-f(x),-f'(x)y)$. First we look for area-preserving reversibilities.

In next theorem we prove that, under the usual conditions considered in relation to integrability around a critical point, the unique possible APR of a classical Li\'enard system is the mirror symmetry w. r. to the $y$ axis. 

\begin{theorem} \label{lierev} 
Let $f,g \in C^\infty \big((a,b),\R \big)$, $a < 0 < b$, with $f(0)=0$, $f'(x) $ and  $xg(x) $ both positive in $(a,b) \setminus \{0\} $. Then the system (\ref{syslie}) has an APR if and only if: 
$$
f(-x) = -f(x), \qquad g(-x) = -g(x).
$$
In such a case, the unique APR is
$$
\sigma(x,y) = (-x,y).
$$
\end{theorem}
{\it Proof.}  Necessity.

 If an APR $\sigma$ exists, then, setting $\sigma(x_1,y_1) = (x_2,y_2) $, by theorem \ref{SRD} one has
$$
\left\{  \matrix{ -f(x_2) = f(x_1)   \hfill  \cr -f'(x_2)y_2 = -f'(x_1)y_1 .}  \right.
$$
Since $f(x)$ is invertible, one has
$$
x_2 = f^{-1}\big( -f(x_1) \big).
$$From the second equality one has
$$
y_2 = \frac{ f'(x_1)}{ f'(x_2)}\, y_1 = \frac{ f'(x_1)}{ f'( f^{-1}\big( -f(x_1) \big))}\, y_1  .
$$
Let us define the functions $\alpha$ and $\beta$ as follows,
\begin{equation} \label{alfabetarev}
\alpha(x) = f^{-1}\big( -f(x) \big), \qquad \beta(x) = \frac{ f'(x)}{ f'(\alpha(x))}  .
\end{equation} 
 From what above, if an APR $\sigma$ exists, it has the form
$$
\sigma(x,y) = \left( \alpha(x),  \beta(x) y  \right)  .
$$
$\sigma$ is an involution because $\alpha$ is an involution,
$$
\alpha(\alpha(x)) =  f^{-1}\left(-f \left( f^{-1}\left( -f(x) \right)  \right)\right) = x,
$$
 and
$$
\beta(\alpha(x) ) \beta(x) y = \frac{ f'(\alpha(x))}{ f'(\alpha(\alpha(x)))}  \frac{ f'(x)}{ f'(\alpha(x))}  y = 
 \frac{ f'(\alpha(x))}{ f'(x)}  \frac{ f'(x)}{ f'(\alpha(x))}  y = y.
 $$
The jacobian matrix of $\sigma$ is
$$
 J_\sigma = \left( \matrix{
 \alpha'(x) & 0 \hfill  \cr \beta'(x) y & \beta(x)
 }  \right).
$$
If $\sigma$ is an APR, then $|\det J_\sigma| = 1$:
$$
1 =  \left| \alpha'(x)  \beta(x) \right|=  \left| \bigg( f^{-1}\big( -f(x) \big) \bigg)'\frac{ f'(x)}{ f'(\alpha(x))} \right| =   
 $$  $$
 = \left|  \frac {-f'(x)}{f'(f^{-1}(-f(x)))} \frac{ f'(x)}{ f'(\alpha(x))} \right| =  \left(   \frac{ f'(x)}{ f'(\alpha(x))}   \right)^2 .
 $$
Together with the sign condition on $f'(x)$, this implies 
$$
f'(\alpha(x)) =  f'(x).
$$
As a consequence, one has
$$
\beta(x) =  \frac{ f'(x)}{ f'(\alpha(x))} = 1.
$$
From the equality $  \alpha'(x)  \beta(x) = \pm 1 $ one has $ \alpha'(x)  = \pm 1$, hence $\alpha( x ) =\pm x + c_\alpha$, $c_\alpha \in \R$. 
Moreover, 
$$
c_\alpha = \alpha(0) = f^{-1}\big( -f(0) \big) = 0, 
$$
hence $\alpha( x ) = \pm x $. The choice $\alpha(x) = x$ leads to the identity map, which is not an APR. Choosing   $\alpha(x) = -x$ leads to
$$
-x  = \alpha(x) =  f^{-1}\big( -f(x) \big) \iff  f(-x) = - f(x).
$$
From what above the reversibility has the form $\sigma(x,y) = \left( -x,  y  \right) $.

\medskip

Sufficiency.

Assume $\sigma(x,y) =(-x,y)$, and let us check the reversibility condition (\ref{Jrev}). Setting $F(x,y) = (y,-g(x) - yf(x))$, one has
$$
F(\sigma(x,y)) =F((-x,y)) = (  y  , -g(-x ) - y f(-x )  ) ,
$$
Then
$$
- J_\sigma(x,y) \cdot F(x,y) = (  y,  g(x) + y  f(x)  ).
$$
The condition (\ref{Jrev}) holds if and only if $g(-x)= - g(x)$ and $f(-x)= - f(x)$.

\hfill  $\clubsuit$ \\

If we look for an APS under the hypotheses of theorem \ref{lierev}, we are lead to the following system,
$$
\left\{  \matrix{ -f(x_2) =- f(x_1)   \hfill  \cr -f'(x_2)y_2 = -f'(x_1)y_1 .}  \right.
$$
Since $f(x)$ is assumed to be strictly increasing in $(a,b)$, the first equation leads to $x_2=x_1$, then the second one to $y_2 = y_1$, hence the unique APS is the identity.
As a consequence, in order to find non-trivial  area-preserving symmetries, we need to change hypotheses.

\medskip

By corollary \ref{notinv}, if an APS $\sigma$ exists, then the map $\Delta(x,y) = (D(x,y),  \dot D(x,y))=$ 
$(-f(x), - f'(x)y)$ is not  locally invertible at a $\sigma$-fixed point. Since its jacobian matrix has the form
$$
 \left( \matrix{
 -f'(x) & 0 \hfill  \cr -f''(x) y &- f'(x)
 }  \right),
$$
the jacobian determinant is $f'(x)^2$ and has to vanish at every $\sigma$-fixed point.
As usual, we assume $O = (0,0)$ to be the unique critical point of (\ref{syslie}) in some strip $(a,b) \times \R$. Since $\sigma$ takes fixed points into fixed points, one has $\sigma(O) =O$. This requires $f'(0) = 0$. In order to simplify the problem, we assume $f(x)$ to be strictly decreasing in $(a,0)$ and strictly increasing in $(0,b)$. Such assumptions were already considered in the study of Li\'enard systems, for instance when studying limit cycle uniqueness.

\begin{theorem} \label{liesym} 
Let $f,g \in C^\infty \big((a,b),\R \big)$, $a < 0 < b$, $f'(x) < 0$ in $(a,0)$, $f'(x) > 0$ in $(0,b)$, $xg(x) >0$ in $(a,b) \setminus \{0\} $. Then the system (\ref{syslie}) has a non-trivial APS $\sigma$ if and only if: 
$$
f(-x) = f(x), \qquad g(-x) = -g(x).
$$
In such a case, the unique APS is
$$
\sigma(x,y) = (-x,-y).
$$

\end{theorem}
{\it Proof.} Necessity.

 By possibly replacing $a$ or $b$ with a point closer to 0, we may assume to have $f\big((a,0)\big) = f\big((0,b)\big) $.
Let us set  $Z:= f\big((a,0)\big) = f\big((0,b)\big) $. 
By hypothesis, for every $l \in Z$ there exist exactly two points $x_l^-$ and $x_l^+$ such that $f(x_l^-) = f(x_l^+) = l$.
Since $f(x)$ is invertible both in $(a,0]$ and in $[0,b)$, we denote by $f_-^{-1}$,  resp. $f_+^{-1}$, the inverse functions of the restrictions of $f(x)$ to $(a,0]$, resp. $[0,b)$.

Let us define the function $\alpha: (a,b) \to (a,b)$ as follows,
\begin{equation} \label{alfasym}
\alpha(x) = \left\{  \matrix{f_+^{-1}\big(  f_-(x) \big), \quad x \in (a,0]  \hfill  \cr  f_-^{-1}\big(  f_+(x) \big), \quad x \in [0,b)   }  \right.
\end{equation} 
Such a definition is well-posed also for $x=0$, since $f(0)=f_-(0)=f_+(0)$. 
It is easy to show that $\alpha$ is an involution of $(a,b)$ onto itself, taking $(a,0)$  into $(0,b)$ and vice-versa. For $x\in (a,0)$, one has
$$
\alpha'(x) = \bigg( f_+^{-1}\big(  f_-(x) \big) \bigg)' = f_+^{-1}{}'\big(  f_-(x) \big)f_-'(x) = \frac {f_-'(x)} { f_+'\left( f_+^{-1} \big(  f_-(x) \big) \right)}
=  \frac {f_-'(x)} { f_+'\left( \alpha(x) \right)} .
$$
Similarly, for $x\in (0,b)$, one has
$$
\alpha'(x) = \frac {f_+'(x)} { f_-'\left( \alpha(x) \right)} .
$$
In both cases one can write
$$
\alpha'(x) = \frac {f'(x)} { f'\left( \alpha(x) \right)} < 0,
$$
since if $x\in(a,0)$ then $\alpha(x) \in (0,b)$ and  vice-versa.

\medskip

 As in theorem \ref{lierev}, one has $\Delta(x,y) = (D(x,y),\dot D(x,y)) = (-f(x),-f'(x)y)$.
If an APS $\sigma$ exists, then, setting $\sigma(x_1,y_1) = (x_2,y_2) $, by theorem \ref{SRD} one has
$$
\left\{  \matrix{ -f(x_2) =- f(x_1)   \hfill  \cr -f'(x_2)y_2 = -f'(x_1)y_1 .}  \right.
$$
If $\sigma$ is not the identity, the first equality implies 
$$
x_2 = \alpha (x_1).
$$
From the second equality one has
$$
y_2 = \frac{ f'(x_1)}{ f'(x_2)}\, y_1  = \frac{ f'(x_1)}{ f'( \alpha (x_1))}\, y_1  .
$$
Let us define the function $\beta$ as in theorem \ref{lierev},
\begin{equation} \label{betasym}
 \beta(x) = \frac{ f'(x)}{ f'(\alpha(x))}  .
\end{equation} 
 From what above, if an APS $\sigma$ exists, it has the form
$$
\sigma(x,y) = \left( \alpha(x),  \beta(x) y  \right)  .
$$
$\sigma$ is an involution because 
$$
\beta(\alpha(x) ) \beta(x) y = \frac{ f'(\alpha(x))}{ f'(\alpha(\alpha(x)))}  \frac{ f'(x)}{ f'(\alpha(x))}  y = 
 \frac{ f'(\alpha(x))}{ f'(x)}  \frac{ f'(x)}{ f'(\alpha(x))}  y = y.
 $$
The jacobian matrix of $\sigma$ is
$$
 J_\sigma = \left( \matrix{
 \alpha'(x) & 0 \hfill  \cr \beta'(x) y & \beta(x)
 }  \right).
$$
If $\sigma$ is an  APS, then $|\det J_\sigma| = 1$:
$$
1 =  \alpha'(x)  \beta(x) =    \left(   \frac{ f'(x)}{ f'(\alpha(x))}   \right)^2 .
 $$
Since  $f'(x)$ and $f'(\alpha(x))$ have opposite signs but at the origin, this implies 
$$
f'(\alpha(x)) = - f'(x).
$$
As a consequence, one has
$$
\beta(x) =  \frac{ f'(x)}{ f'(\alpha(x))} = -1.
$$
From the equality $  \alpha'(x)  \beta(x) = 1 $ one has $ \alpha'(x)  = - 1$, hence $\alpha( x ) =-x + c_\alpha$. 
Moreover, 
$$
c_\alpha = \alpha(0)  = 0, 
$$
hence $\alpha( x ) =-x $. Then, from the first equality in (\ref{alfasym}), one has, for $ x \in (a,0)$,
$$
 - x = f_+^{-1}\big(  f_-(x) \big) \iff  f_+(-x) = f_-(x )\iff  f(-x) = f(x).
 $$
 Similarly for   $ x \in (0,b)$, hence  the symmetry has the form $\sigma(x,y) = \left( -x,  -y  \right) $.

\medskip

Sufficiency.

Assume $\sigma(x,y) =(-x,-y)$, and let us check the symmetry condition (\ref{Jsymm}). Setting $F(x,y) = (y,-g(x) - yf(x))$, one has
$$
F(\sigma(x,y)) =F((-x,-y)) = (  -y  , -g(-x ) + y f(-x )  ) .
$$
Then
$$
 J_\sigma(x,y) \cdot F(x,y) = ( - y,  g(x) + y  f(x)  ).
$$
The condition (\ref{Jsymm}) holds if and only if $g(-x)=-g(x)$ and $f(-x)=f(x)$.

\hfill  $\clubsuit$ \\

The condition $f'(x) < 0$ in $(a,0)$, $f'(x) > 0$ in $(0,b)$ can be replaced by exchanging the inequalities, $f'(x) > 0$ in $(a,0)$, $f'(x) < 0$ in $(0,b)$, without changes in the proof. In fact, if $x(t)$ is a solution to the equation
$$
\ddot x + f(x) \dot x + g(x) = 0,
$$
then $x(-t)$   is a solution to the equation
$$
\ddot x - f(x) \dot x + g(x) = 0.
$$
The existence of a symmetry for the first equation is equivalent to the existence of a symmetry for the second equation.

 %%%%%%%%%%%%%%%%%%%%%%%%%%%%%%%%%%%%%% 

\enddocument